\documentclass[10pt,leqno]{article}
\usepackage{amsmath}
\usepackage{amsfonts}
\usepackage{amsthm}
\usepackage{indentfirst}

\newcommand{\N}{\mathbb{N}}

\newcommand{\C}{\mathbb{C}}
\newcommand{\f}{\rightarrow}
\newcommand{\deb}{\bar\partial}
\newcommand{\de}{\partial}
\newcommand{\K}{K\"{a}hler}
\newcommand{\lmb}{\lambda}
\newcommand{\ol}{\operatorname{Hol}}
\newcommand{\hilb}{\mathcal{H}}
\newcommand{\ep}{\varepsilon}
\newcommand{\Aut}{\operatorname{Aut}}
\newcommand{\Isom}{\operatorname{Isom}}
\newcommand{\D}{\mathcal{D}}

\newtheorem{theorem}{\indent\rm T h e o r e m\;}[section]

\newtheorem{definition}{\indent\rm D e f i n i t i o n\;}[section]

\renewenvironment{proof}{\indent\rm P r o o f.\;}{\hfill $\square$ \\ \indent}

\makeatletter                                                  %
\renewcommand*{\@seccntformat}[1]{
  \csname the#1\endcsname\;-                                   %
}                                                              %
\renewcommand{\section}{\@startsection{section}{1}{0mm}        %
   {1.5\baselineskip}
   {1\baselineskip}
   {\indent\normalfont\normalsize\bfseries}
   }                                                           %
\renewcommand*{\@seccntformat}[1]{
  \normalfont\bfseries\csname the#1\endcsname\;-               %
}                                                              %
\renewcommand\subsection{\@startsection                        %
  {subsection}{2}{0mm}
  {1.5\baselineskip}
  {1\baselineskip}
  {\indent\normalfont\normalsize\itshape}}
\renewcommand*{\@seccntformat}[1]{
  \normalfont\bfseries\csname the#1\endcsname\;-               %
}                                                              %
\renewcommand\subsubsection{\@startsection                     %
  {subsubsection}{2}{0mm}
  {1.5\baselineskip}
  {1\baselineskip}
  {\indent\normalfont\normalsize\texttt}}
\makeatother                                                   %
\begin{document}
\thispagestyle{empty}



\centerline{\sc\large{\textbf{A bounded homogeneous domain and }}}
 \centerline{\sc\large{\textbf{a projective manifold are not relatives}}}
\vspace {1.5cm}

\begin{center}
{\sc\large Roberto Mossa}
\end{center}

\renewcommand{\thefootnote}{\fnsymbol{footnote}}

\footnotetext{
Research partially supported by GNSAGA (INdAM) and MIUR of Italy}

\renewcommand{\thefootnote}{\arabic{footnote}}
\setcounter{footnote}{0}

\vspace{1,5cm}
\begin{center}
\begin{minipage}[t]{10cm}
\small{
\noindent \textbf{Abstract.}
Let $M_1$ and $M_2$ be two \K\ manifolds. Following \cite{dilo} one says that 
 $M_1$ and $M_2$ are 
{\em relatives} if they share a non-trivial \K\ submanifold
$S$, namely, if there exist two holomorphic and isometric immersions
(\K\ immersions) $h_1: S\rightarrow M_1$ and $h_2: S\rightarrow
M_2$. Our main results in this paper  is  Theorem \ref{mainteor} 
where we  show that a bounded homogeneous domain with a homogeneous
\K\ metric and a  projective \K\  manifold
(i.e. a projective manifold  endowed with the restriction of the  Fubini--Study metric)
are not relatives.
Our result is a generalization of the result  obtained in \cite{dilo}
for the Bergman metrics.
\medskip

\noindent \textbf{Keywords.} Diastasis, homogeneous bounded domains, \K\ metric, relatives.
\medskip

\noindent \textbf{Mathematics~Subject~Classification~(2010):}
53D05, 53C55.

}
\end{minipage}
\end{center}


\bigskip


\section{Introduction}
The study of \K\ immersions (holomorphic and isometric immersions) started with E. Calabi in his seminal paper \cite{calabi} where he gave necessary and sufficient conditions for the existence of a \K\ immersion of a  finite dimensional   \K\ manifold into a complex space form. In  particular  he proved that two complex space forms with curvature  of different
sign cannot be \K\ immersed one into another.  Moreover, he proved that  for complex space forms of  the same type, just projective spaces can be embedded between themselves in a non trivial way
by using Veronese mappings. Almost 25 years later  M. Umehara  \cite{um} proved  that two  finite dimensional complex space forms with holomorphic sectional curvatures of different signs
cannot be relatives. Recall the definition of relatives:
\begin{definition}(\cite[Definition 1.1]{dilo})
Let $r\geq 1$ be an  integer.
Two \K\ manifolds $M_1$ and $M_2$ are said to be {\em $r$-relatives} if they have in common a complex $r$-dimensional \K\ submanifold
$S$, i.e. there exist two \K\  immersions $h_1:S\rightarrow M_1$ and
$h_2:S\rightarrow M_2$.
Otherwise, we say that $M_1$ and $M_2$ are
not relatives. 
\end{definition}
Recently, A. J. Di Scala and A. Loi \cite{dilo} proved the following:
\begin{theorem}\label{thm andrea} (\cite[Theorem 1.2]{dilo}) A bounded domain $D\subset {\C}^n$
endowed with its Bergman metric and a  projective \K\  manifold endowed with the restriction of  the Fubini--Study metric are not relatives.
\end{theorem}
Our main result is the following  theorem which  generalizes  Theorem \ref{thm andrea} when the Bergman metric is homogeneous.
Recall that a $n$-dimensional  bounded homogeneous domain $(\Omega, g)$ is a bounded domain of $\C^n$ endowed with the \K\ metric $g$ such that the group $G=\Aut (\Omega) \cap \Isom  (\Omega,g)$ act transitively on it. Here $\Aut (\Omega)$ denotes the group of invertible holomorphic maps of $\Omega$ and $\Isom (\Omega, g)$  the group of isometries of 
$(\Omega, g)$.

\begin{theorem}\label{mainteor}
A bounded homogeneous domain $\left ( \Omega, g \right)$ and a  projective \K\  manifold endowed with the restriction of the Fubini--Study metric are not relatives.
\end{theorem}
The proof of our this theorem  is based on a recent result  \cite{lm1} obtained by the author jointly  with A. Loi. We point out that our result is of local nature, i.e. no assumptions are used about the compactness or completeness of the manifolds involved.

\section{Proof of Theorem \ref{mainteor}}
In order to prove Theorem \ref{mainteor} it is enough to show that a bounded homogeneous domain and  the complex projective space $\C P^m$ are not relatives.
Let $\omega$ be the \K\ form associated to $g$. 
It is well-know that there exists a globally defined K\"ahler potential $\Phi$ for $g$ i.e. $\omega=\frac{i}{2}\de\deb \Phi$. Indeed, $\Omega$ is pseudoconvex being biholomorphically equivalent to a Siegel domain (see, e.g. \cite{vinberg} for a proof)  and so the existence of a global potential follow by a classical result of Hormander (see \cite{bulletin}) asserting that the equation $\deb u = f$ with $f$ $\deb$-closed form has a global solution on pseudoconvex domains
(see also the proof of Theorem 4 in \cite{dlh}, for an explicit construction of the potential $\Phi$ following  the ideas developed in \cite{dorf}).
Consider the associated weighted Hilbert space $\hilb_{\lmb\Phi}$ of square integrable holomorphic functions on $\Omega$, with weight $e^{-\lmb\Phi}$, namely
\begin{equation}\label{hilbertspacePhi}
\hilb_{\lmb\Phi}=\left\{ f\in\ol(\Omega) \ | \ \int_\Omega e^{-\lmb\Phi}|f|^2\frac{\omega^n}{n!}<\infty\right\},
\end{equation}
In \cite{lm1} it is proven that for $\lmb>0$ large enough $\hilb_{\lmb\Phi}\neq \{0\}$. Fixed such $\lmb$, let $K_{\lmb\Phi}(z, w)$ be its reproducing kernel.
We can define the $\ep$-function:
\begin{equation}\label{epsilon}
\ep_{\lmb g}(z)=e^{-\lmb\Phi(z)}K_{\lmb\Phi}(z, z),
\end{equation}
this function does not depend on the potential $\Phi$, it depends only on the constant $\lmb$ and on the metric $g$ (see \cite{lm1} for details). Moreover, it is invariant with respect to the action of the  Lie  group $G$. Since $G$ acts transitively on $\Omega$, 
it follows that $\ep_{\lmb g}=C$ is constant (for $\lmb$ large enough) and $\log K_{\lmb\Phi} (z,z) $ is a \K\ potential for the metric $\lmb g$. By analytic continuation we have
\begin{equation}\label{epsilon1}
\ep_{\lmb g}(z,w)=e^{-\lmb\Phi(z,w)}K_{\lmb\Phi}(z, w)=C>0,
\end{equation}
 and so $K_{\lmb\Phi}(z, w)$ never vanishes. Then, fixed a point $z_0$, the function
\begin{equation*}
\psi(z,w)=\frac{K_{\lmb\Phi}(z,w)K_{\lmb\Phi}(z_0,z_0) }{K_{\lmb\Phi}(z,z_0)K_{\lmb\Phi}(z_0,w)}.
\end{equation*}
is well defined. Observe that $\psi(z_0,w)=\psi(z,z_0)=1$ and that
\begin{equation*}
\D(z)=\log\psi(z,z)
\end{equation*}
is a globally defined \K\ potential for $g$ (actually $\D(z)$ is  the \emph{diastasis} centered in $z_0$, see Calabi in \cite{calabi} for details and further property about the diastasis function). We can now consider the Hilbert space $\hilb_{\lmb\D}$ 
 given by:
\begin{equation*}
\hilb_{\lmb\D}=\left\{ f\in\ol(\Omega) \ | \ \int_\Omega e^{-\lmb\D}|f|^2\frac{\omega^n}{n!}<\infty\right\},
\end{equation*}

Let us denote $K_{\lmb\D}(z, w)$ its reproducing kernel, as  the $\ep$-function does not depend on the \K\ potential, by \eqref{epsilon1} we have
\begin{equation}\label{epsilon2}
\ep_{\lambda g}(z,w)=e^{-\lmb\D(z,w)}K_{\lmb\D}(z, w)=C,
\end{equation}
where $\D(z,w)$ is the analytic continuation of $\D(z)$.
In particular
\begin{equation*}
K_{\lmb\D}(z_0, w)=K_{\lmb\D}(z, z_0)=C
\end{equation*}
and so $\hilb_{\lmb\D}$ contains the constant functions and by boundedness of $\Omega$ all polynomials belong to $\hilb_{\lmb\D}$. In particular $\hilb_{\lmb\D}$ contains the sequence $\{z_1^k\}_{k \in \N}$, by applying the Gram-Schmidt orthonormalization  procedure we get a sequence $P=\{P_k\}_{k\in \N}$ of orthonormal polynomials in the variable $z_1$.

Consider now the {\em coherent states map} (see \cite{{lm1}}) $\varphi:\Omega\rightarrow\C P^{\infty}$ from $\Omega$ into the infinite dimensional
complex projective space $\C P^{\infty}$ given by
\begin{equation}\label{states}
\varphi: \Omega\f \C P^{\infty},\  z\mapsto [P_0(z_1),P_1(z_1), \dots ,F_0(z), F_1(z), \dots ].
\end{equation}
where $\{P_0(z_1), \dots ,F_0(z), \dots\}$ is an orthonormal basis of $\hilb_{\lmb\D}$ obtained completing  $P$ to an orthonormal basis.  Since $K_{\lmb\D}(z, z)= \sum_{k=0}^\infty |P_k(z_1)|^2 + |F_k(z)|^2$ by \eqref{epsilon2} we see that the map \eqref{states}   is well-defined.
Moreover, the constancy of $\ep$ also implies that 
$\varphi^*g_{FS}=\lambda g$,
where $g_{FS}$ is the Fubini--Study  metric on $\C P^{\infty}$ (see \cite{rawnsley} for a proof).
In  other words,  the metric $\lambda g$ is projectively induced via
the  coherent states map.

Assume now,  by contradiction, that $(\Omega, g)$ is $r$-relative (for some positive integer $r\geq 1$) to the complex projective space $\C P^m$
(the \K\ metric on $M$ is induced by the Fubini-Study metric through the immersion $j$).
Then  we can assume $r=1$ and that there exists an open subset $S$ of $\C$ through the origin and two \K\ immersions $f:S \f D$ and $h:S \f \C P^m$.

Consider the \K\ map $\varphi \circ f:S \f \C P^\infty$ 
\begin{equation*}
\varphi \circ f \left(\xi\right)=\left[P_0\left(f_1\left(\xi\right)\right),\,P_1\left(f_1\left(\xi\right)\right), \dots ,\,F_0\left(f\left(\xi\right)\right),\, F_1\left(f\left(\xi\right)\right), \dots \right]
\end{equation*}
were $f_1$ is the first component of $f$. Without loss of generality assume $\frac{\de f_1}{\de \xi}(0)\neq 0$.
We claim  that $P_0\left(f_1\left(\cdot\right)\right),P_1\left(f_1\left(\cdot\right)\right),\dots$ are linearly independent functions on $S$. Let $a_0,\dots,a_q$ be complex numbers such that 
\begin{equation*}
a_0 P_{k_0}\left(f_1\left(\xi\right)\right)+\dots+a_qP_{k_q}\left(f_1\left(\xi\right)\right)=0.
\end{equation*}
By the assumption on $f_1$ it follows that $f_1(S)$ contains an open set of $\C$, therefore $a_0 P_{k_0}\left(z\right)+\dots+a_qP_{k_q}\left(z\right)=0$ for every $z\in \C$. Since $P_{k_0}, \dots, P_{k_q}$ are linearly independent all the $a_k$ must be zero, proving our claim.

 On the other hand, if   $i: \C P^m \f \C P^\infty$ is the standard  totally geodesic embedding, then
 $i \circ h:S \f \C P^\infty$ is a \K\ immersion.
Thus the smallest subspace containing  $i \circ h (S)$ is infinite dimensional, while $\varphi \circ f$  is contained in a $m$-dimensional subspace, this is in contrast  with the Calabi's rigidity theorem (see \cite{calabi}), as wished.
\endproof



\vspace{0.5cm} \indent {\it
A\,c\,k\,n\,o\,w\,l\,e\,d\,g\,m\,e\,n\,t\,s.}\; I wish to thank Prof. Andrea Loi for his interest in my work and various stimulating discussions.

\bigskip
\begin{center}

\end{center}

\bigskip
\bigskip
\begin{minipage}[t]{10cm}
\begin{flushleft}
\small{
\textsc{Roberto Mossa}
\\*Laboratoire de Mathématiques Jean Leray (UMR 6629) CNRS,
\\*2 rue de la Houssinière
(B.P. 92208)
\\*44322 Nantes Cedex 3, France
\\*e-mail: roberto.mossa@gmail.com
}
\end{flushleft}
\end{minipage}

\end{document}